\documentclass[runningheads]{llncs}
\input{001_setup}

\begin{document}

\title{The wrought iron beauty of Poncelet loci\vspace{-1em}}
\author{Dan Reznik\orcidID{0000-0002-0542-6014}}
\authorrunning{D. Reznik}
\institute{Data Science Consulting, Rio de Janeiro, Brazil\\
\email{dreznik@gmail.com}}
\maketitle

\vspace{-0.75cm}

\begin{abstract}
We've built a web-based tool for the real-time interaction with loci of Poncelet triangle families. Our initial goals were to facilitate exploratory detection of geometric properties of such families. During frequent walks in my neighborhood, it appeared to me Poncelet loci shared a palette of motifs with those found in wrought iron gates at the entrance of many a residential building. As a result, I started to look at Poncelet loci aesthetically, a kind of generative art. Features were gradually added to the tool with the sole purpose of beautifying the output. Hundreds of interesting loci were subsequently collected into an online ``gallery'', with some further enhanced by a graphic designer. We will tour some of these byproducts here. An interesting question is if Poncelet loci could serve as the basis for future metalwork and/or architectural designs.
\end{abstract}
\vspace{-0.75cm}

\setlength{\epigraphwidth}{0.95\textwidth}
\epigraph{The mathematician’s patterns, like the painter’s or the poet’s, must be beautiful; the ideas, like the colors or the words, must fit together in a harmonious way. Beauty is the first test: there is no permanent place in the world for ugly mathematics. --G.H. Hardy, ``A mathematician’s 
apology''}

\vspace{-0.75cm}
\section{Introduction}
The front doors of apartment buildings lining many a street in Rio de Janeiro, Brazil, are bedecked with wrought iron\footnote{Iron was popular for ornamental applications in the 19th and early 20th century, given its high ductility and structural capabilities. It has now been by and large replaced by steel.} masterpieces. Their curvaceous designs seem to strike a perfect balance between beauty and sturdiness.

\begin{figure}
    \centering
    \includegraphics[width=\textwidth]{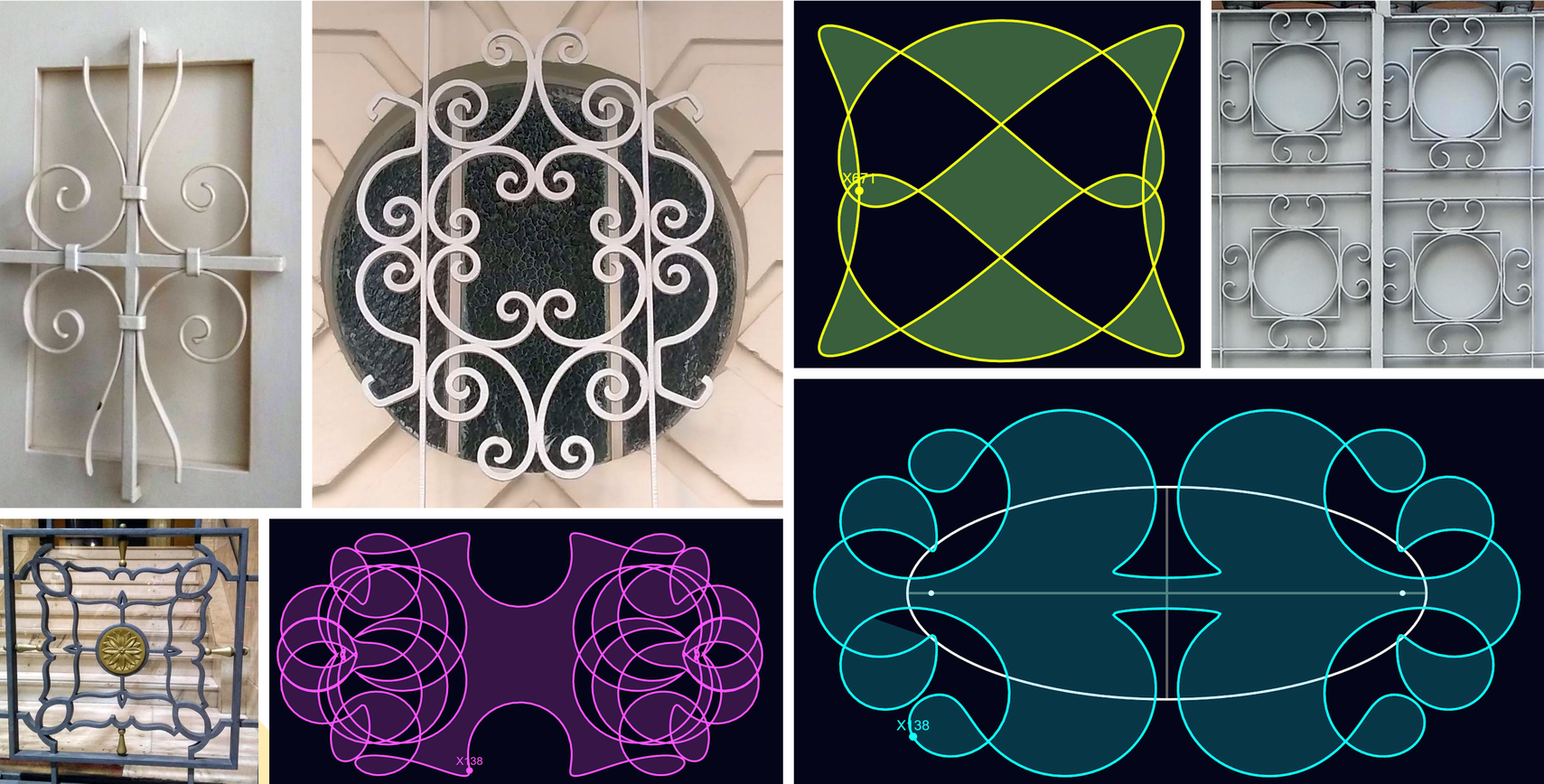}
    \caption{Wrought iron gates and loci of Poncelet triangles, hinting at common design motifs.}
    \label{fig:pastiche}
\end{figure}

As one experiments with the curves swept by points attached to Poncelet triangle families (more details in \cref{sec:poncelet}), one becomes aware of certain common motifs linking these two worlds, e.g.,   four-fold symmetry, harmonious spirals, tiling, recursion, etc.,
see  \cref{fig:pastiche}.

\begin{figure}
    \centering
    \includegraphics[width=\textwidth]{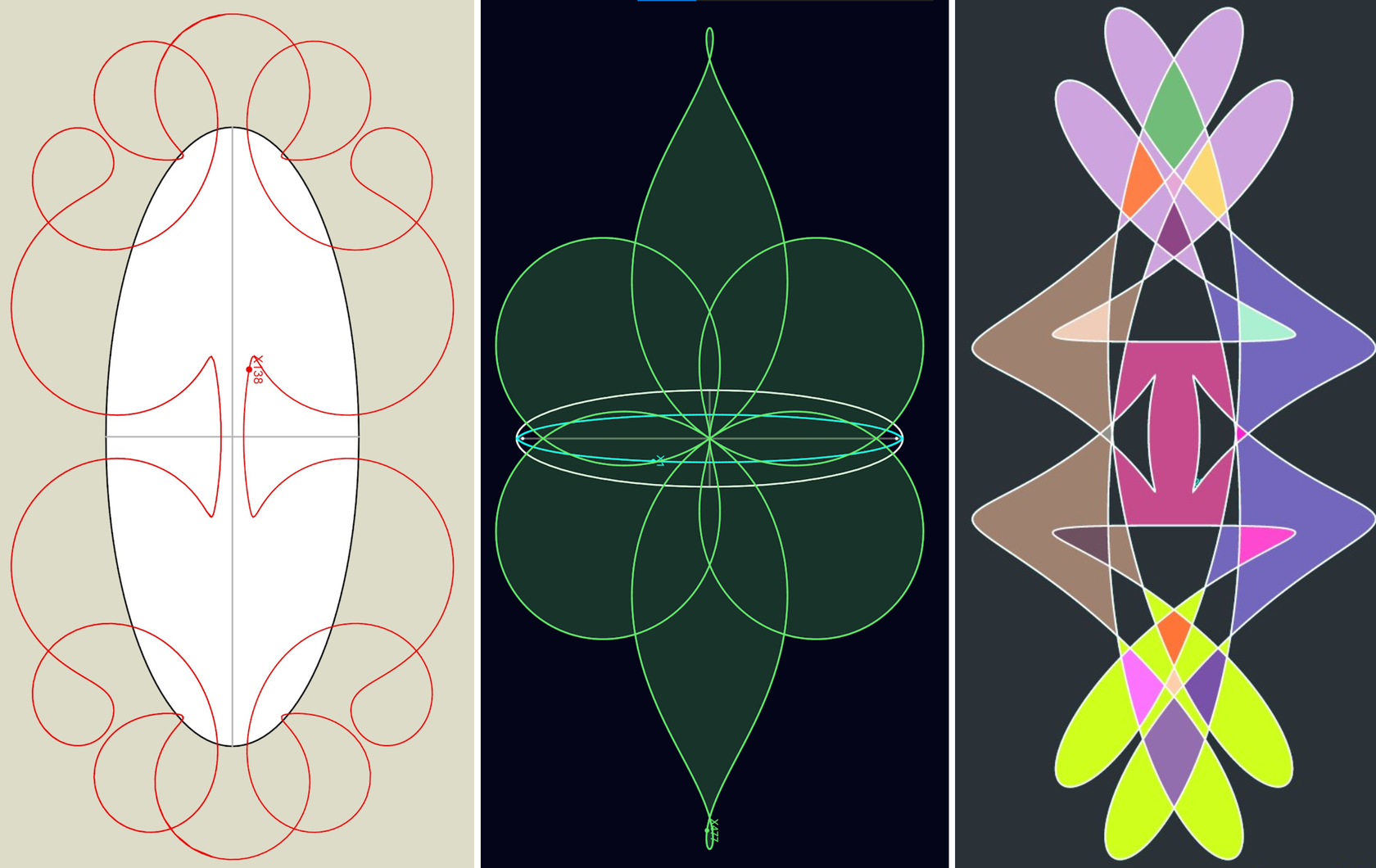}
    \caption{Beautifying loci: from wireframe (left), to a thick curve against a dark background (middle), to a region-colored design (right). Experimental parameters used to produce them can be shared as a URL and/or exported as vector graphics for further processing by a graphic artist.}
    \label{fig:evolution}
\end{figure}

We had built a web-based tool to interact with Poncelet triangles and their loci, see \cite{darlan2021-app,reznik2021-locus-app-tutorial}. As new experiments were set up, we would often stumble upon a new ornate locus. We quickly recognized in Poncelet loci a kind of aesthetic talent. Gradually, we added features geared at beautifying, coloring, and sharing such curves, see \cref{fig:evolution}. At the same time, we started to collect hundreds of aesthetically-pleasing finds into slideshow and video ``galleries'' \cite{reznik2021-artful-long,reznik2021-artful-short,reznik2021-artful-youtube}. Leafing through these gives one an idea of the wide design palette possible with Poncelet loci.

Our goal here is to explore the ingredients involved in generating such loci. In the next sections we (i) review Poncelet's porism, (ii) review the basic geometry of a triangle's notable points (whose loci we will sweep), (iii) explore some features of our locus-rendering app, and (iv) tour a few the artful byproducts of Poncelet loci.

\subsubsection*{Related work}
The fields of kinetic art \cite{popper2003} and computer-generated (generative) art \cite{dreher2015-art} have been evolving for several decades. Works that explore the connection between classical art and geometry include  \cite{gutruf2010-vermeer,karger98-classical,lordick2021-vermeer}. The usage of dynamic geometry tools for mathematical discovery is beautifully expounded in  \cite{sharp2015-artzt}. 
In \cite{odehnal2011-poristic,pamfilos2004,skutin2013,zaslavsky2001-poncelet,zaslavsky2003-trajectories} loci of triangle centers are studied over various 1d triangle families, Poncelet or otherwise. Works \cite{reznik2020-ballet,reznik2020-intelligencer} follow in their footsteps and identify new curious properties and invariants of Poncelet families. Proofs that loci of certain triangle centers in over the confocal family are ellipses appear in \cite{corentin2021-circum,garcia2019-incenter,garcia2020-new-properties,olga14}. A theory of locus ellipticity is slowly emerging, see \cite{helman2021-power-loci,helman2021-theory}. In  \cite{garcia2020-similarity-I,garcia2020-family-ties,reznik2020-similarityII}, similar geometric properties and invariants are used to cluster Poncelet triangles families. Proofs of  experimentally-detected invariants of Poncelet families appear in  \cite{akopyan2020-invariants,bialy2020-invariants,caliz2020-area-product,stachel2021-billiards}.

\section{Loci of Poncelet Triangles}
\label{sec:poncelet}
\begin{figure}
    \centering
    \includegraphics[width=.8\textwidth]{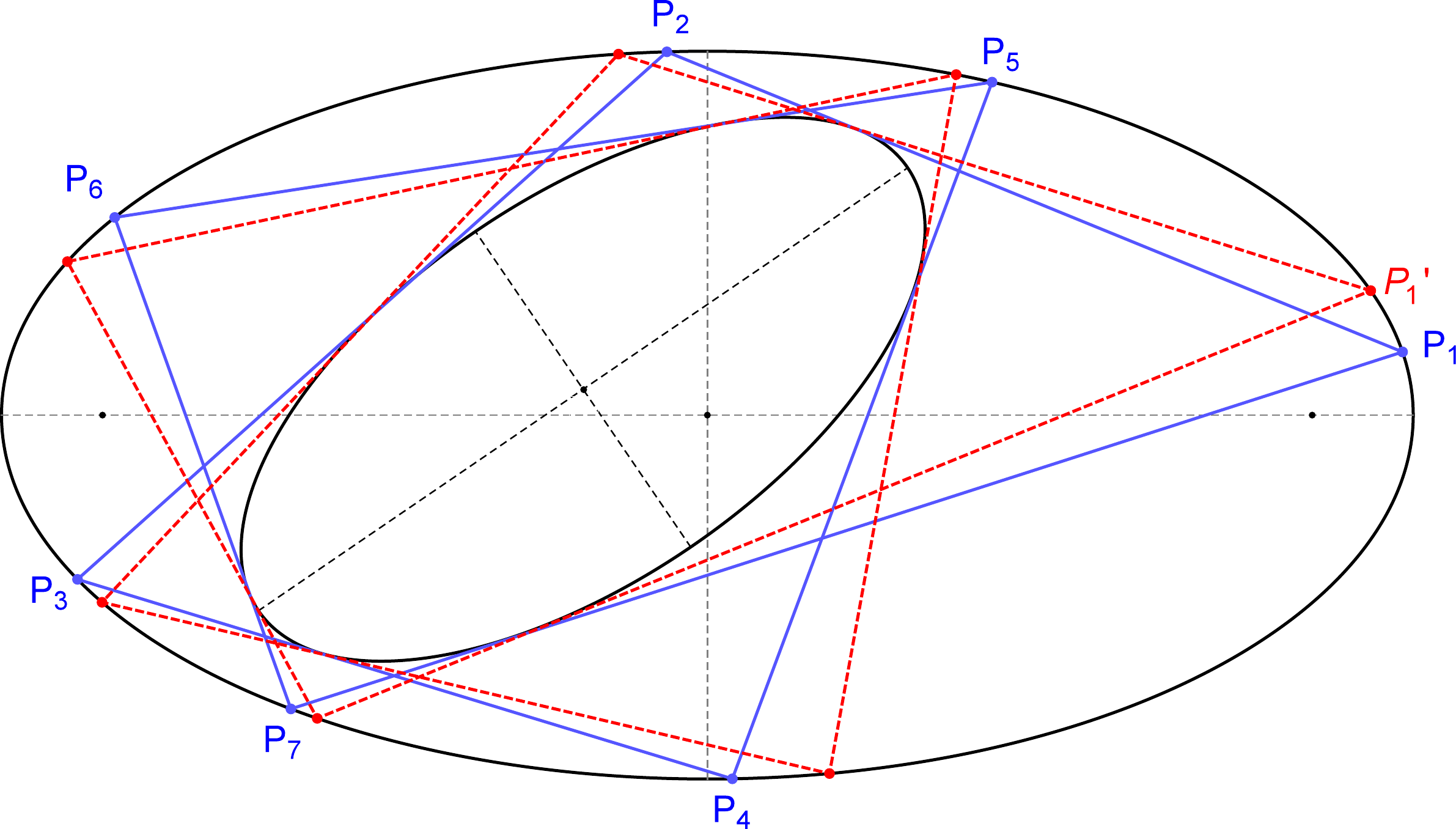}
    \caption{A heptagon $P_1,...,P_7$ (blue) is shown inscribed in an outer ellipse and circumscribing an inner one. In such a case, Poncelet's theorem guarantees that any point of the outer ellipse can be used to construct a 7-sided polygon $P_i'$ similarly inscribed/circumscribed about the two conics. A smooth traversal of the entire family can be seen \href{https://youtu.be/kzxf7ZgJ5Hw}{here}.}
    \label{fig:poncelet}
\end{figure}

Poncelet's closure theorem is illustrated in \cref{fig:poncelet}. It states that given two real conics\footnote{Recall these comprise ellipses, hyperbolas, parabolas, as well as some degenerate specimens \cite{stachel2019-conics}.} $\C,\C'$, if one can draw a polygon with all vertices on $\C$ and with all sides tangent to $\C'$, then a one-dimensional family of such polygons exists. \cite{centina15,dragovic11,bos1987}. 

\begin{figure}
    \centering
    \includegraphics[width=\textwidth]{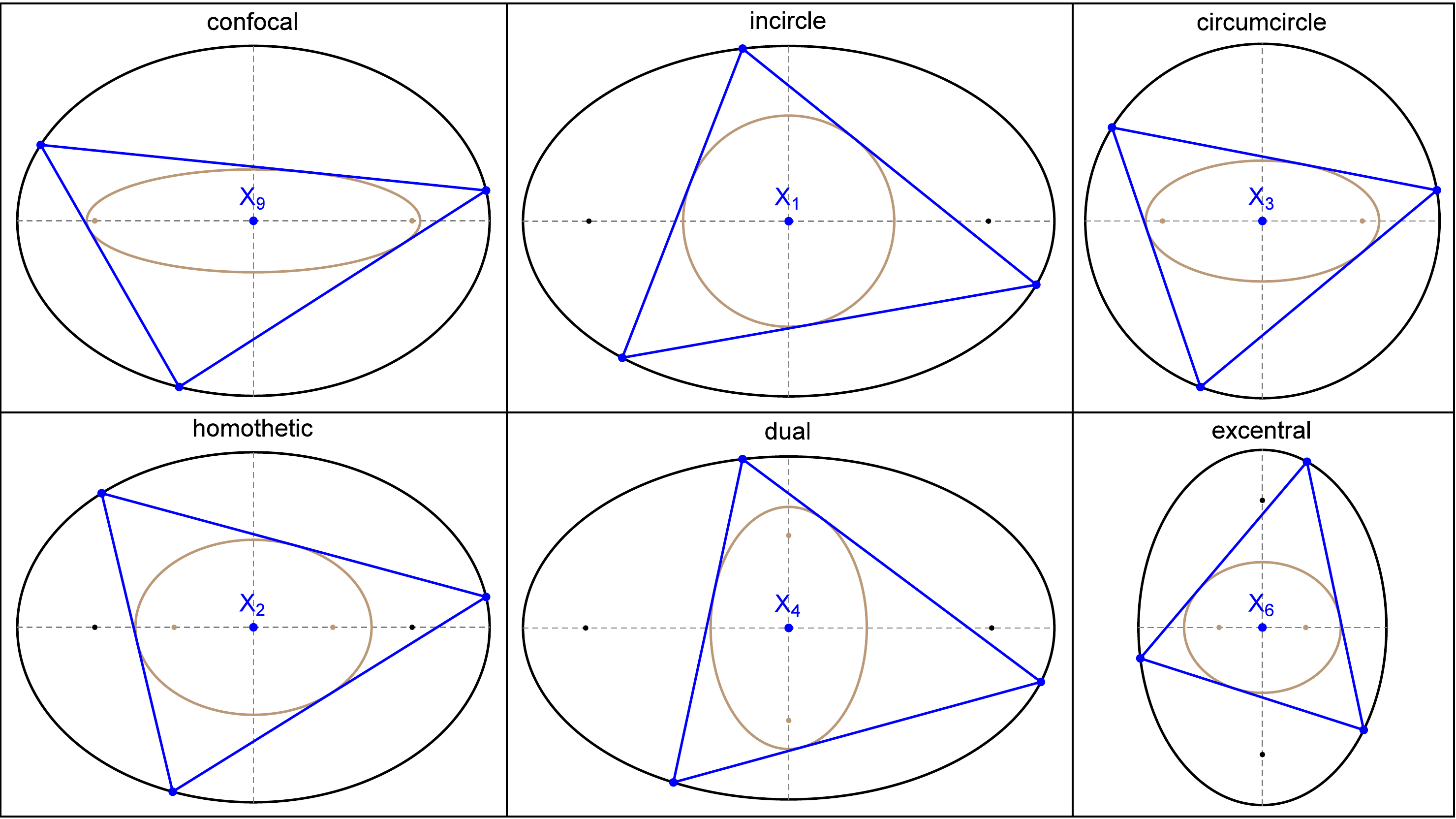}
    \caption{Six Poncelet triangle families inscribed in an ellipse and circumscribing a concentric, axis-parallel inellipse or caustic. We call these ``confocal'', ``incircle'', ``circumcircle'', ``homothetic'', ``dual'', ``excentral''. In each case, a certain triangle center remains stationary at the common center, namely, $X_k$, $k=$9,1,3,2,4,6, respectively. \href{https://youtu.be/14TQ5WlZxUw}{Video}}
    \label{fig:six-caps}
\end{figure}

We will herein focus on families of Poncelet {\em triangles} families. In \cref{fig:six-caps}, six examples are shown of such families interscribed\footnote{This is shorthand for ``inscribed while simultaneously circumscribing''.} between two concentric, axis-aligned ellipses\footnote{In general, the pair of Poncelet conics need not be concentric nor axis aligned.}.

\begin{figure}
\centering
\begin{subfigure}[t]{0.48\textwidth}
 \centering
 \includegraphics[width=\textwidth]{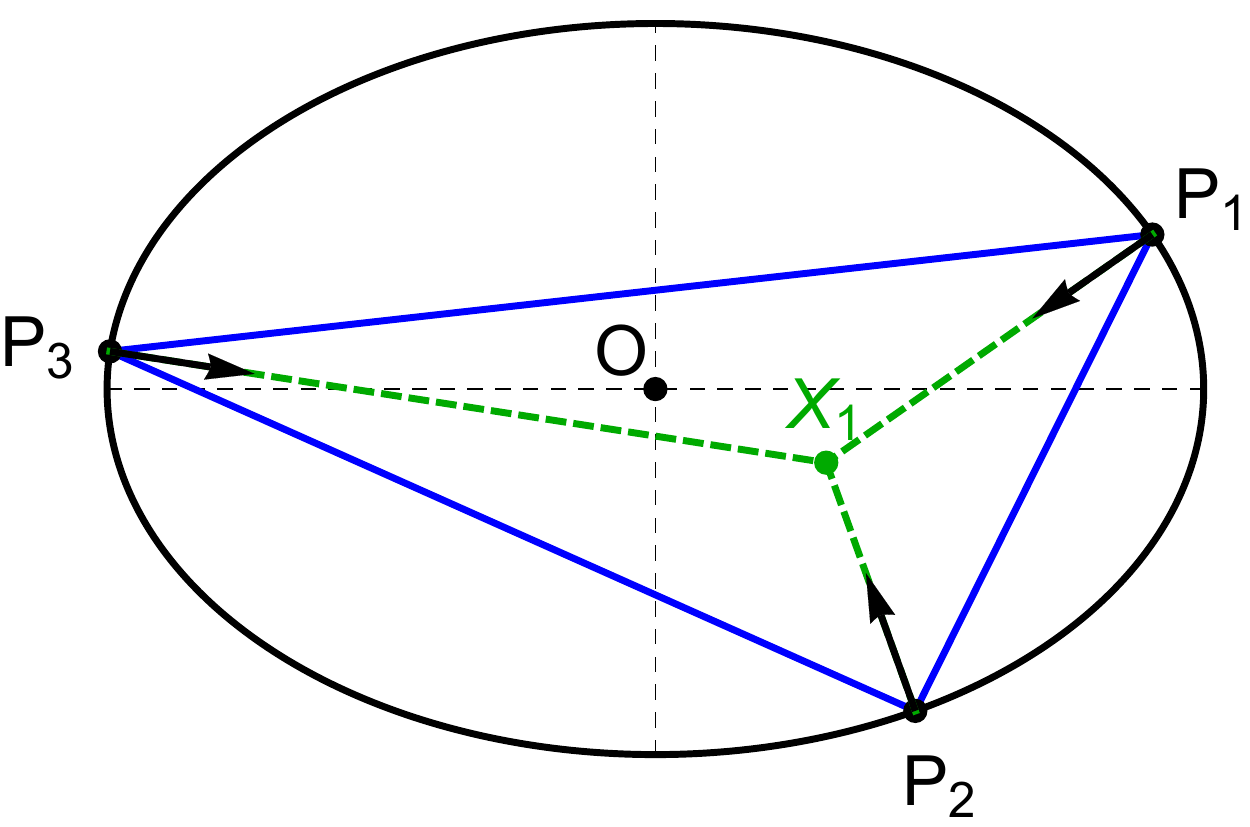}
\end{subfigure}
\hfill
\begin{subfigure}[t]{0.48\textwidth}
 \centering
  \includegraphics[width=\textwidth]{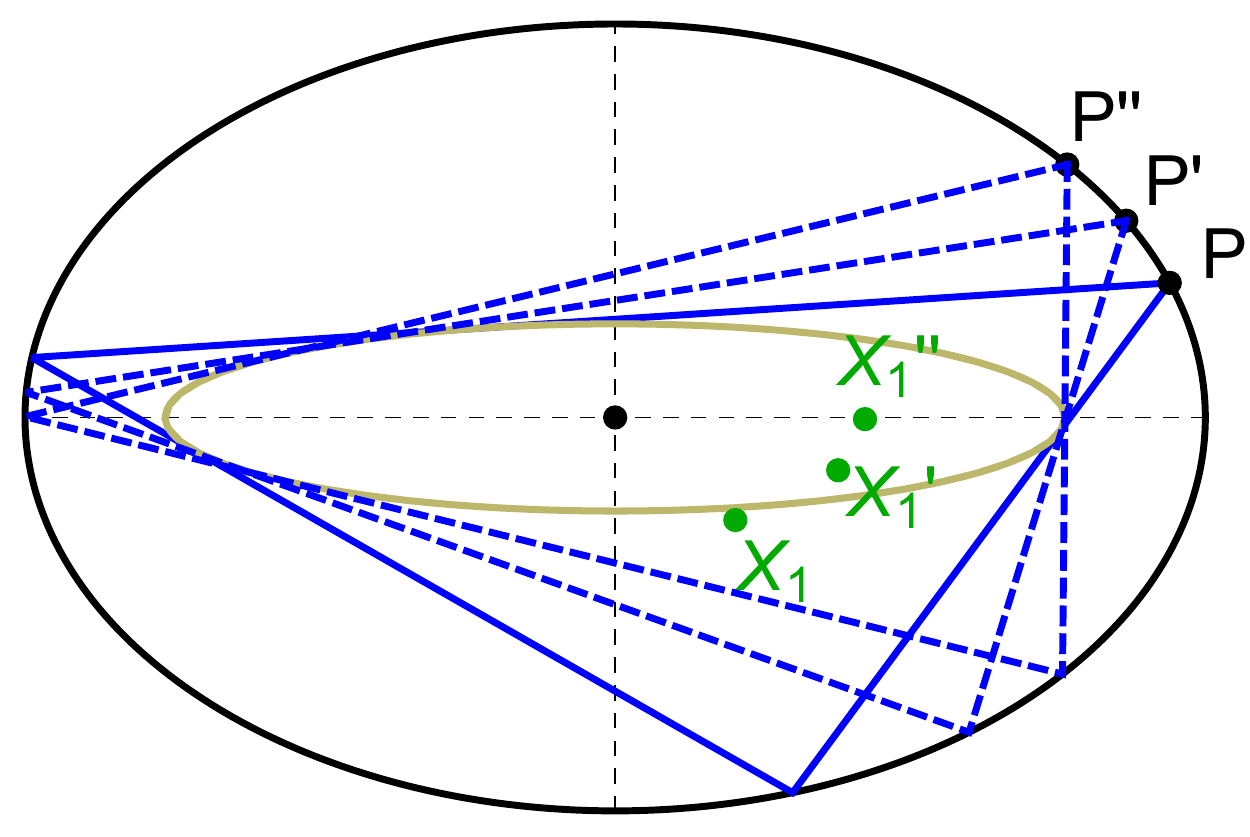}
\end{subfigure}
     \caption{\textbf{Left:} The incenter $X_1$ of a triangle inscribed in an ellipse and circumscribing a confocal one (shown on the right). \textbf{Right}: As one sweeps the triangle family, the incenter changes position: $X_1,X_1',X_1''$, sweeping a locus.  
     \href{https://youtu.be/BBsyM7RnswA}{Video 1}, 
\href{https://youtu.be/Y3q35DObfZU}{Video 2} }
\label{fig:x1}
\end{figure}

\begin{figure}
    \centering
    \includegraphics[trim=0 180 0 20,clip,width=\textwidth,frame]{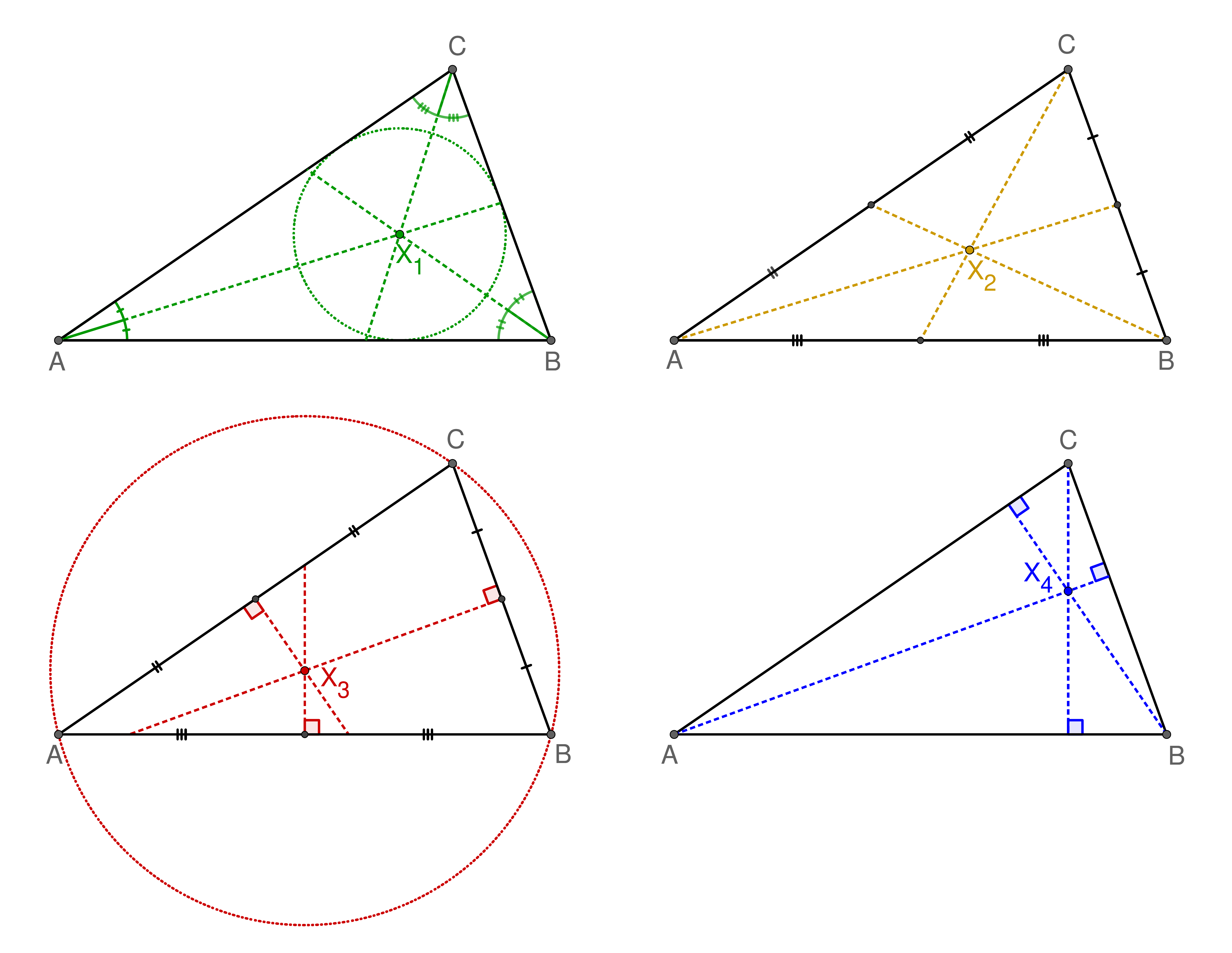}
    \caption{The four notable points of a triangle: (i) top left: the incenter $X_1$, where angular bisectors meet, also the center of the inscribed circle; (ii) top right: the barycenter $X_2$, where medians meet; bottom right: (iii) the circumcenter $X_3$, where perpendicular bisectors meets, also the center of the circumscribed circle; (iv) bottom right: the orthocenter $X_4$, where altitudes meet.}
    \label{fig:notables}
\end{figure}

 Referring to \cref{fig:x1}, a first natural question is: over some particular triangle family, what are curves swept by a {\em notable point}? Recall the four classical notable points of a triangle, shown in \cref{fig:notables}, namely, (i) the incenter $X_1$, (ii) the barycenter $X_2$, (iii) the circumcenter $X_3$, and (iv) the orthocenter $X_4$. The $X_k$ notation conforms with \cite{etc}, where thousands of such points, known as {\em triangle centers}, are specified.

Referring to \cref{fig:notables}(left), an early observation was that over the confocal family, the loci of the four notable points $X_1,X_2,X_3,X_4$ are ellipses. Formal proofs appeared in \cite{olga14,corentin2021-circum,garcia2019-incenter}. Interestingly, other centers can sweep quartics, self-intersecting curves, segments, and be stationary points, see \cref{fig:notables}(right), and \cite{garcia2020-new-properties}.

\begin{figure}
    \centering
    \begin{subfigure}[c]{0.48\textwidth}
    \includegraphics[width=\textwidth]{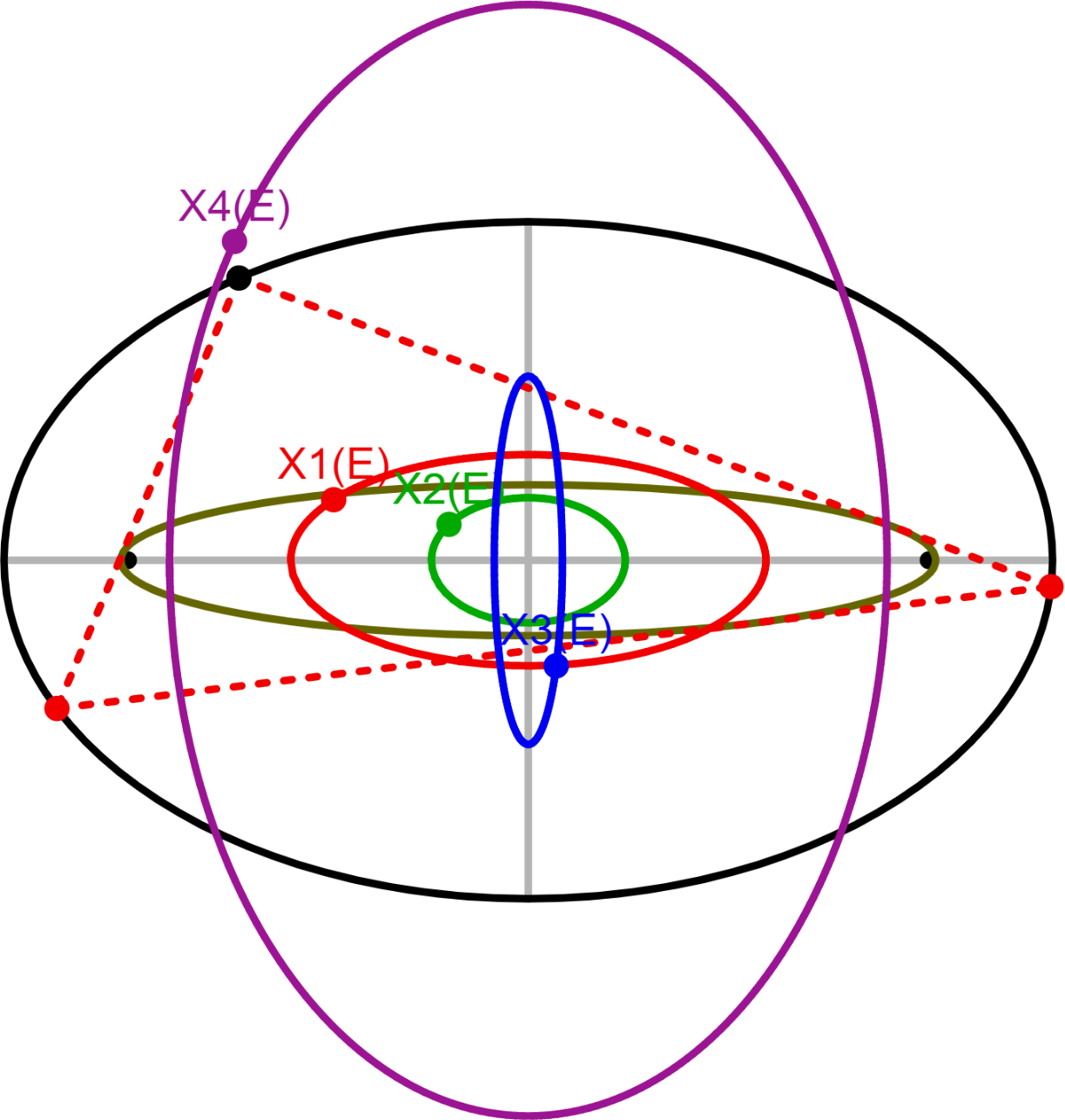}
    \end{subfigure}
    \rulesep
     \begin{subfigure}[c]{0.48\textwidth}
    \includegraphics[width=\textwidth]{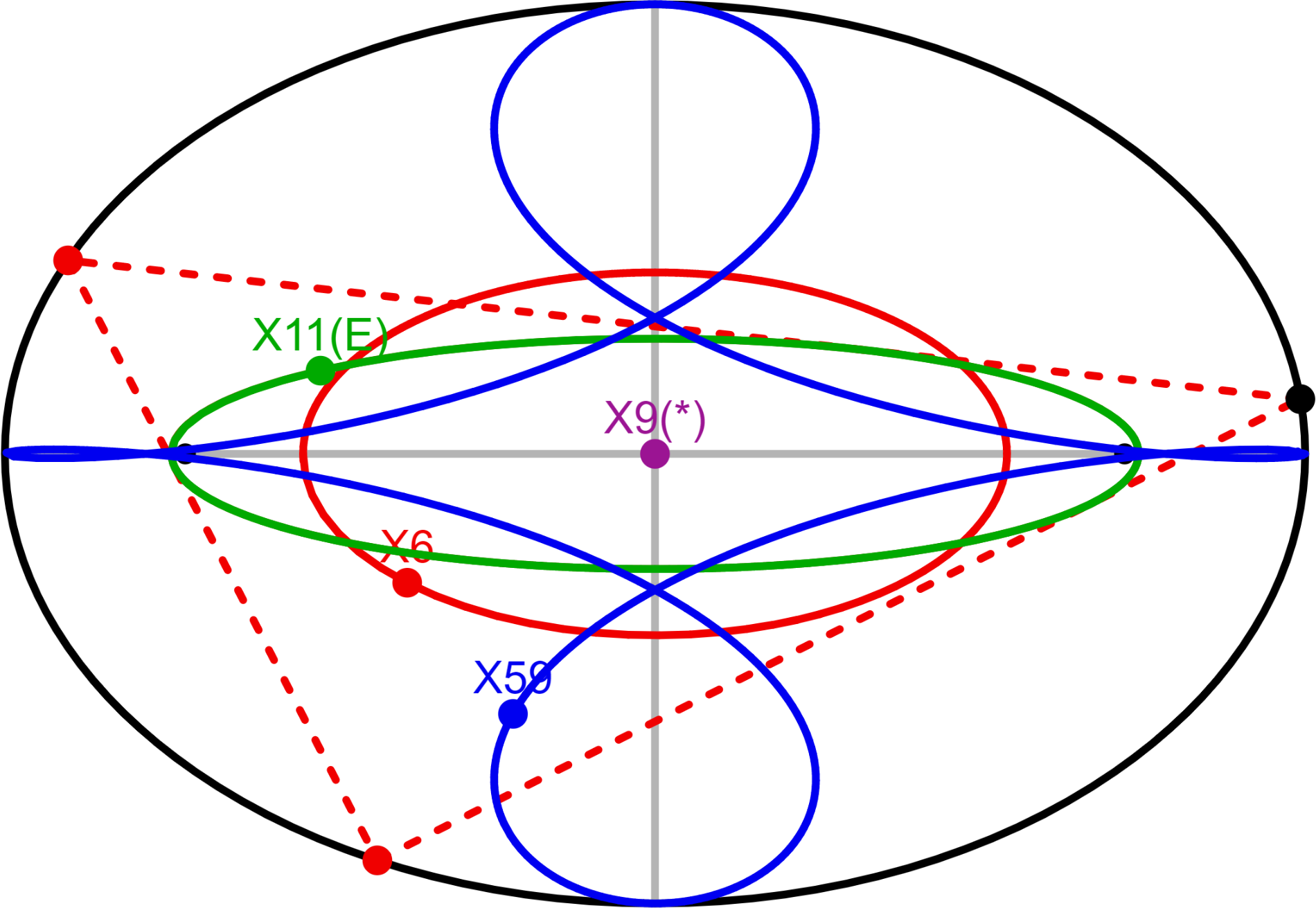}
    \end{subfigure}
    \caption{\textbf{Left:} elliptic loci of the incenter $X_1$ (red), barycenter $X_2$ (green), circumcenter $X_3$, and orthocenter $X_4$ over Poncelet triangles between two confocal ellipses. \textbf{Right:} Over the same family, the mittenpunkt $X_9$ (purple) is  stationary \cite{reznik2020-intelligencer}, the symmedian point $X_6$ (red) sweeps a quartic \cite{garcia2020-new-properties},  the Feuerbach point $X_{11}$ (green) is a curve identical to the caustic (inner ellipse), and $X_{59}$ sweeps a self-intersected curve (blue), studied in \cite{reznik2020-ballet}. \href{https://bit.ly/3qZrhVH}{app}}
    \label{fig:x1234}
\end{figure}

\section{The Locus Rendering App}
\label{sec:app}
\begin{figure}
    \centering
    \includegraphics[width=\textwidth]{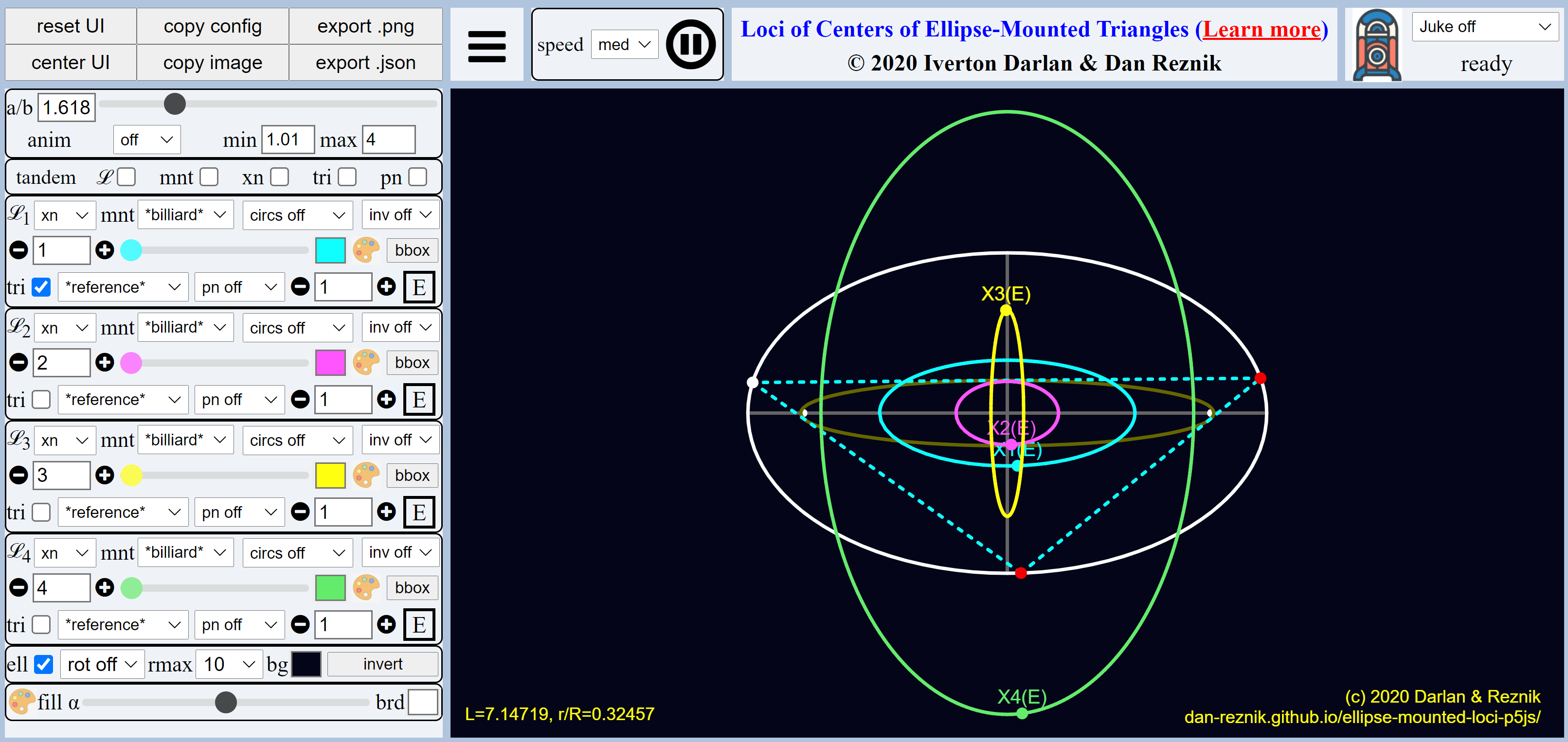}
    \caption{Locus Visualization app to explore 3-periodic families. Shown are the loci of $X_k$, $k=$1,2,3,4, over billiard 3-periodics. The ``(E)'' suffix indicated they are numerically ellipses. \href{https://bit.ly/3yV8caF}{Live}; Also see our tutorial \href{https://bit.ly/3iCKHxn}{playlist}.}
    \label{fig:splash}
\end{figure}

\begin{figure}
    \centering
    \includegraphics[width=\textwidth]{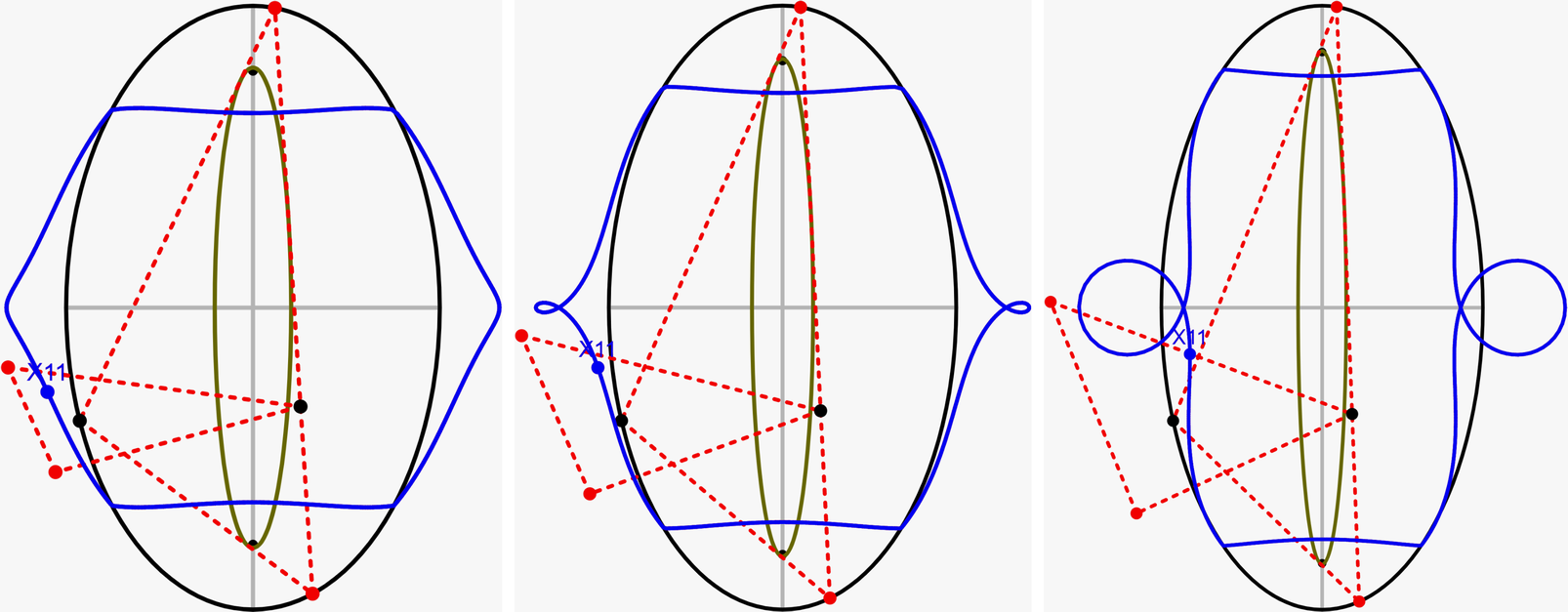}
    \caption{Observing transitions in the topology of tghe locus of the Feuerbach point $X_11$ of the orthic triangle (smaller dashed red) as the aspect ratio of the outer ellipse in the confocal pair is smoothly altered; \href{https://bit.ly/2Zg1caK}{live}.}
    \label{fig:orthic}
\end{figure}

Referring to \cref{fig:splash}, the default view of our tool is the elliptic locus of the incenter over the confocal family. The app can render loci of the first 1000 triangle centers in \cite{etc}, over one dozen Poncelet families, including the ones of \cref{fig:poncelet}. The user can interactively change parameters of the simulation (e.g., aspect ratio of Poncelet ellipses, triangle center tracked, derived triangle being used, etc.), and observe topological changes in the loci being studied. As an example see \cref{fig:orthic}.

\begin{figure}
    \centering
    \includegraphics[width=\textwidth,frame]{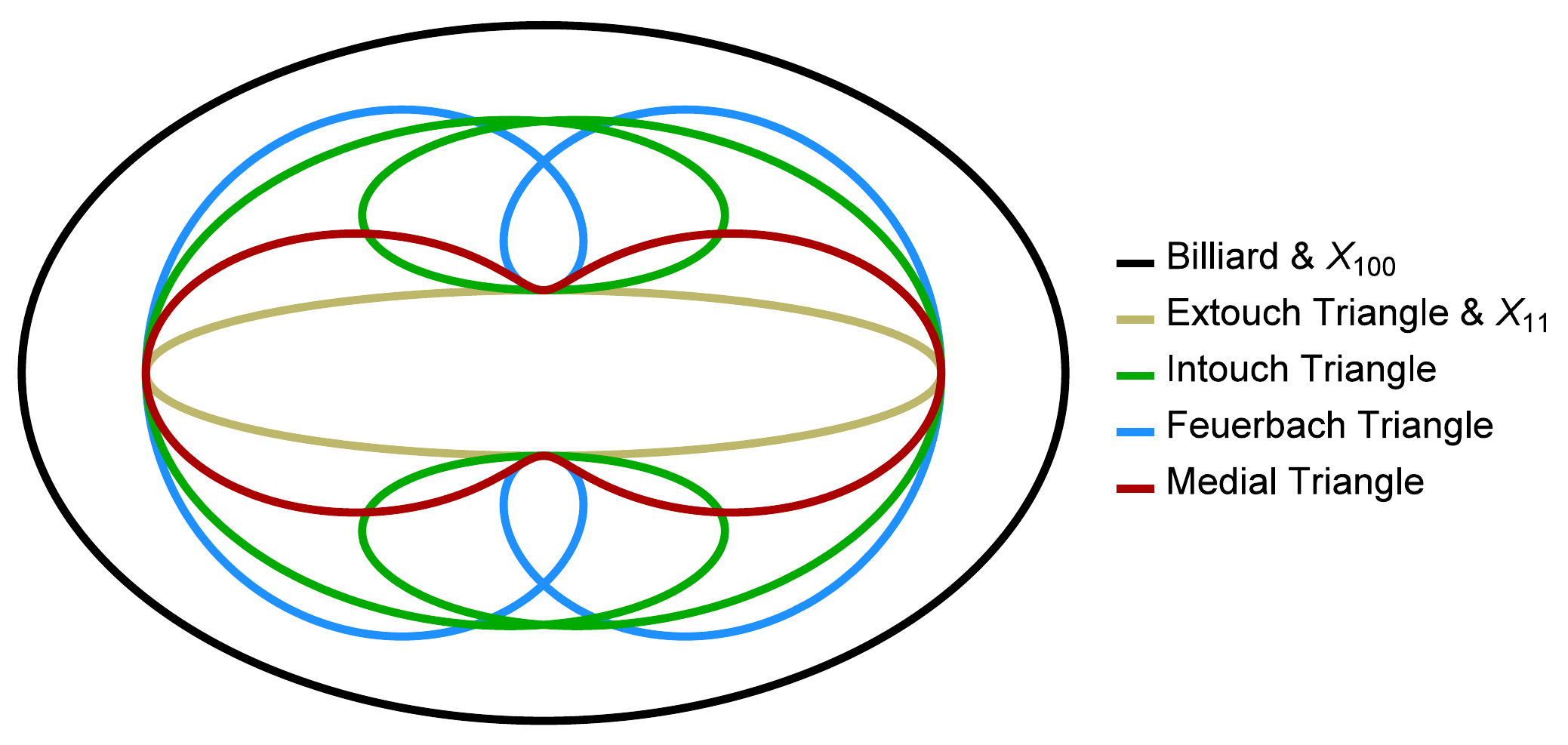}
    \caption{Loci of vertices of well-known derived triangles over the confocal family.  Their construction is defined in \cite{mw}.}
    \label{fig:vertices}
\end{figure}

Loci of triangle centers of {\em derived} triangles can be studied as well, e.g., the orthic, medial excentral, triangles\footnote{These correspond to triangles with vertices which are the feet of altitudes, medians, and where external bisectors meet, respectively.}, etc., see \cite{mw,etc}. As shown in \cref{fig:vertices}, the locus of vertices of the main or derived triangles can be studied as well, further expanding the palette of obtainable curves.

\section{Aestheticizing Poncelet}
\label{sec:art}
\epigraph{As a mathematical discipline travels far from its empirical source, or still more, if it is a second and third generation only indirectly inspired from ideas coming from reality, it is beset with very grave dangers. It becomes more and more purely aestheticizing, more and more purely l'art pour l'art. -- John von Neumann, ``the Mathematician''}

\begin{figure}
    \centering
    \includegraphics[width=\textwidth]{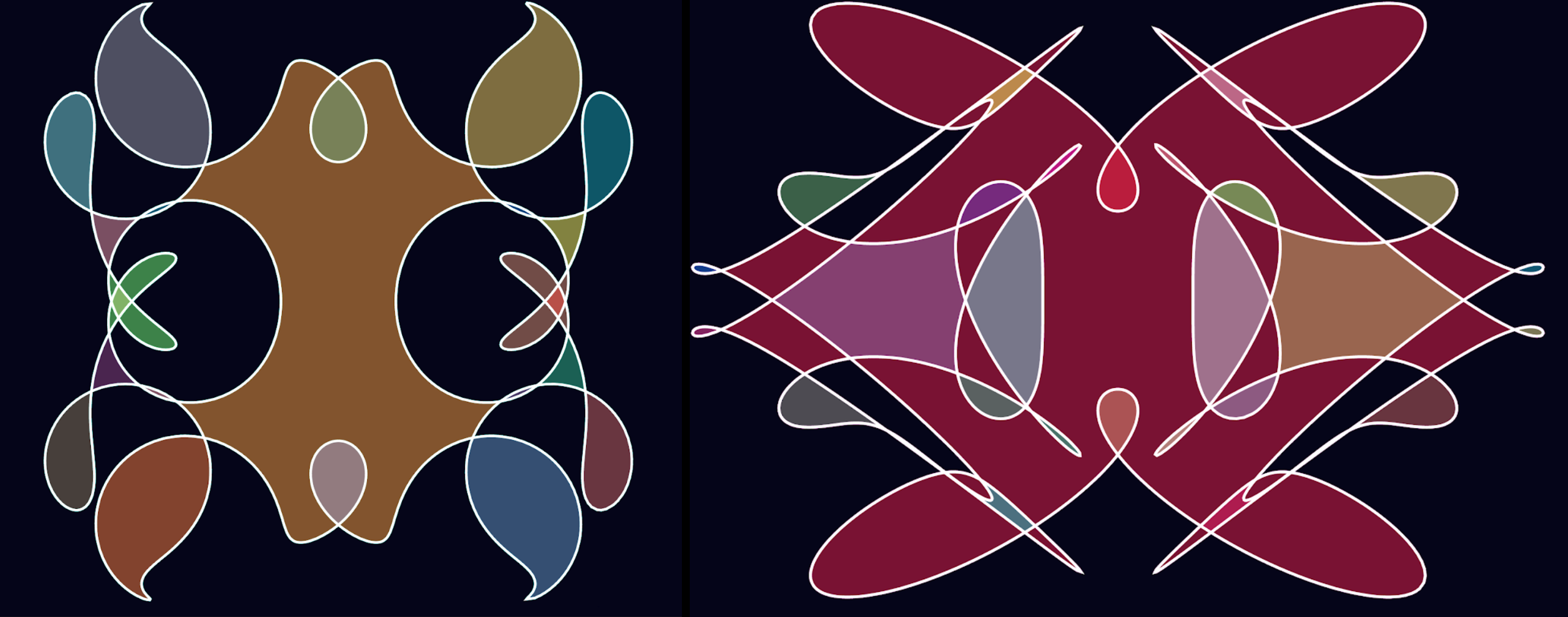}
    \caption{Two ``leafy'' Poncelet loci.}
    \label{fig:leafy}
\end{figure}

\begin{figure}
    \centering
    \includegraphics[width=\textwidth]{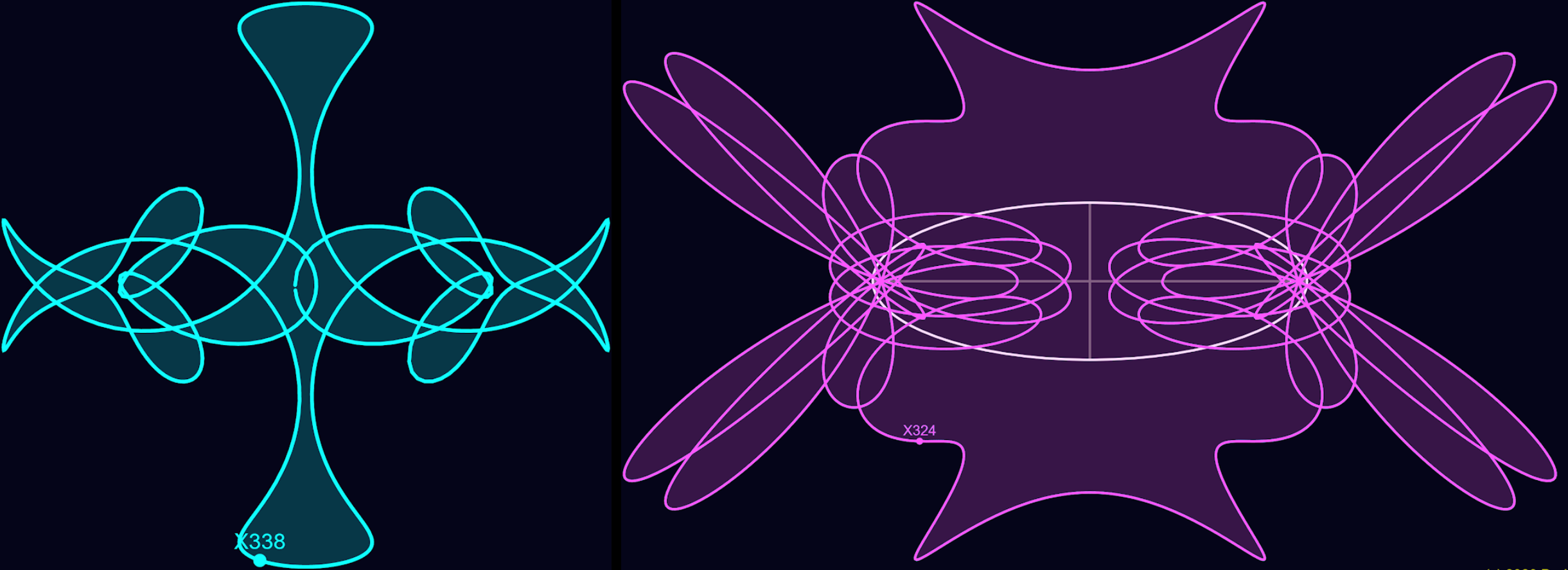}
    \caption{Dragonfly-like loci of Poncelet triangles.}
    \label{fig:libelule}
\end{figure}

\begin{figure}
    \centering
    \includegraphics[trim=0 0 75 0,clip,width=\textwidth]{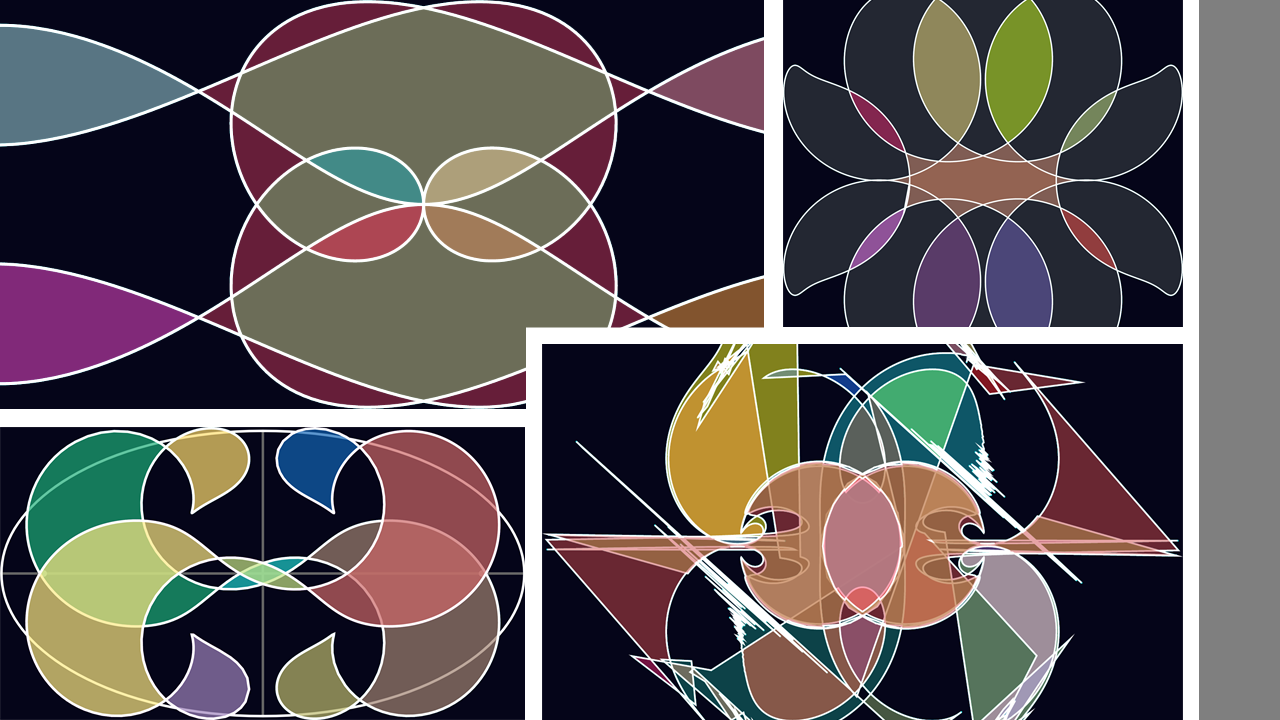}
    \caption{A collage of color-filled loci which can be produced with our web-based tool. For the complete list, see \cite{reznik2021-artful-short,reznik2021-artful-long,reznik2021-artful-youtube}.}
    \label{fig:artful}
\end{figure}

Not exactly following von Neumann's advice, we added features to our locus tool to beautify loci: dark backgrounds, thick outlines, and automatic color-filling of connected regions with a random palette of pastel colors, see \cref{fig:artful,fig:leafy,fig:libelule}. As mentioned above, 100s of such designs are now collected in  \cite{reznik2021-artful-long,reznik2021-artful-short}.

\subsubsection*{Interaction with a digital designer} 
\begin{figure}
    \centering
    \includegraphics[width=\textwidth]{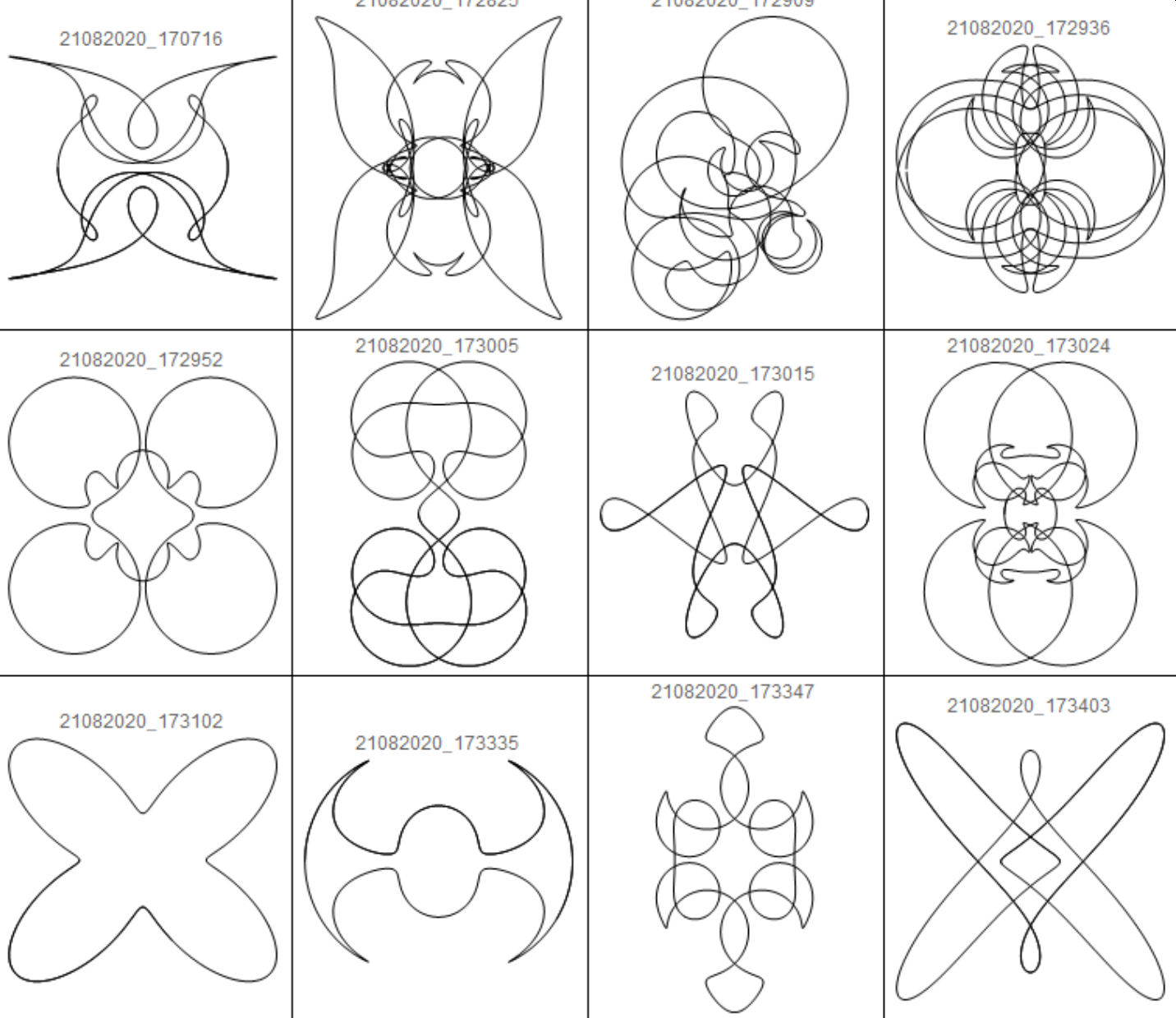}
    \caption{Raw loci sent to digital designer for further coloring and manipulation.}
    \label{fig:raw-loci}
\end{figure}

\begin{figure}
    \centering
    \includegraphics[width=\textwidth]{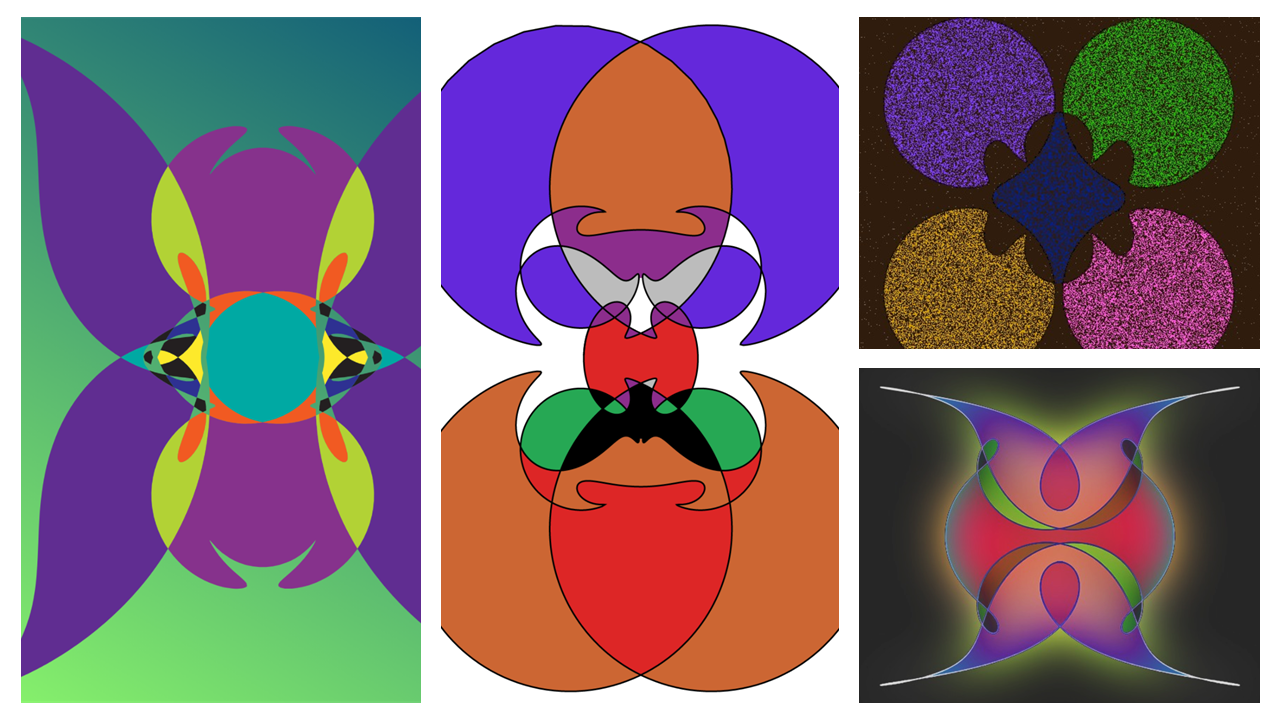}
    \caption{Samples of Regina Reznik's artwork \cite{regina2021-private}.}
    \label{fig:regina}
\end{figure}

Some outlines of interesting loci are shown in \cref{fig:raw-loci}. These were sent in vector format to a graphic designer (my sister \cite{regina2021-private}). Using digital image-editing tools and much creativity, a sample of her work appears in \cref{fig:regina}. 

\section{Conclusion}
\begin{figure}
     \centering
     \begin{subfigure}[c]{0.48\textwidth}
         \centering
         \includegraphics[trim=0 100 0 0,clip,width=\textwidth]{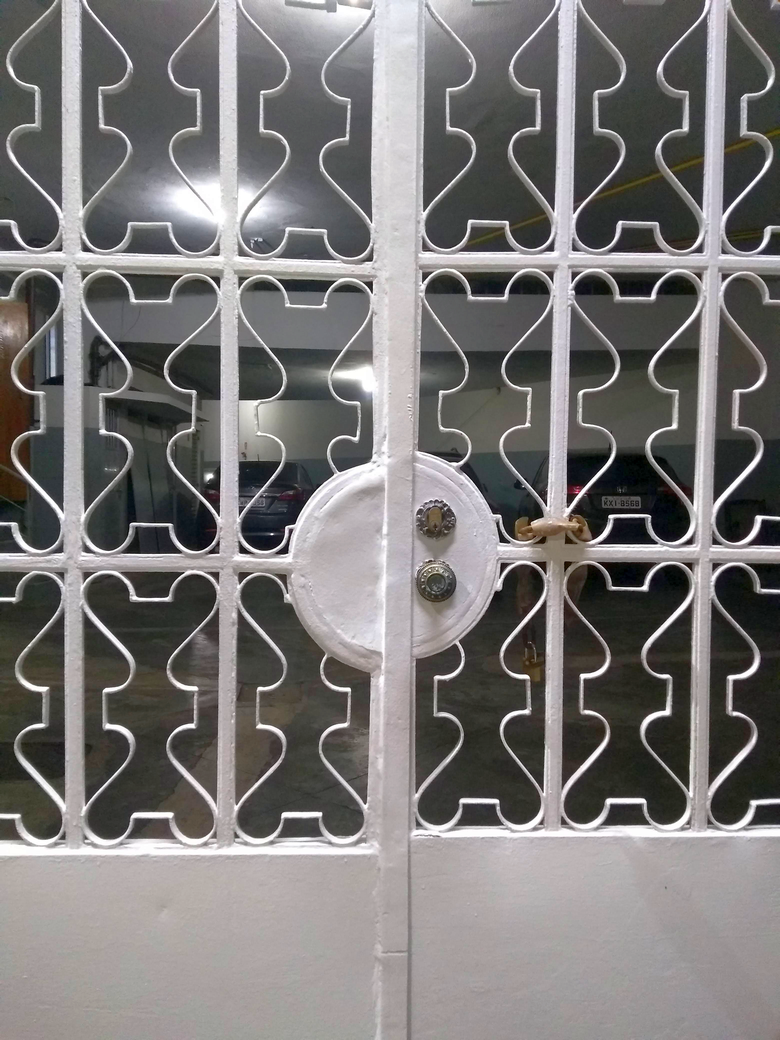}
     \end{subfigure}
     \hfill
     \begin{subfigure}[c]{0.48\textwidth}
         \centering
         \includegraphics[width=\textwidth]{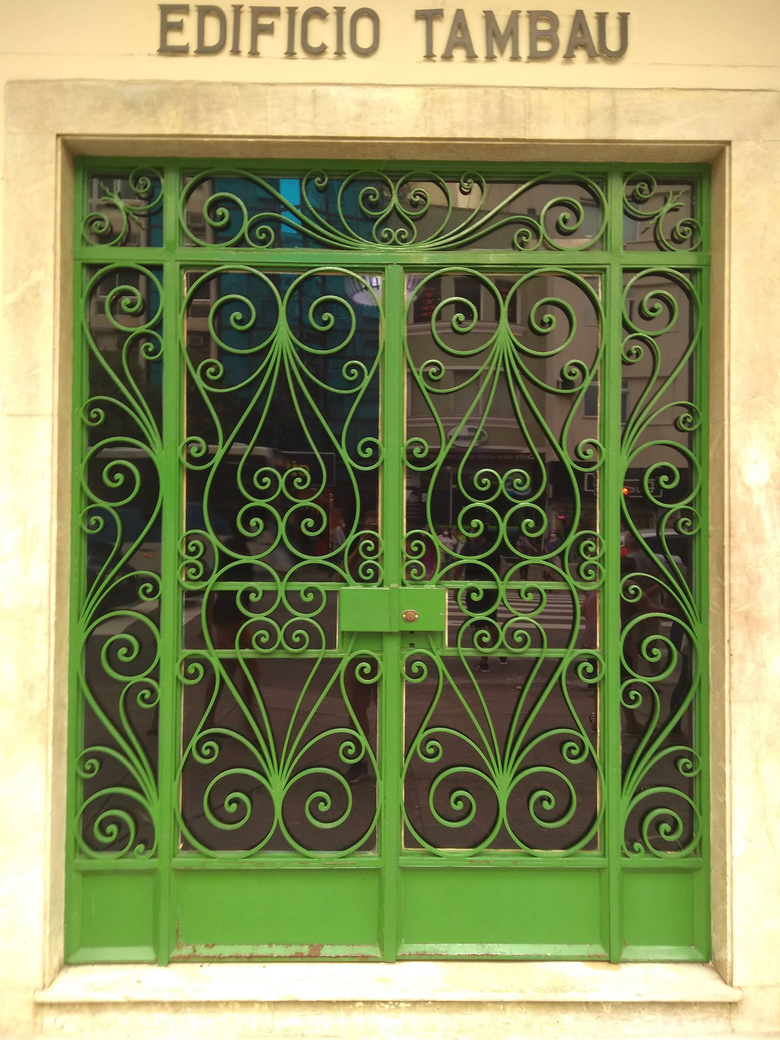}
         \end{subfigure}
        \caption{More wrought iron masterpieces from the streets of Flamengo (left) and Copacabana (right).}
        \label{fig:more-gates}
\end{figure}

Art, architecture, and design have been inspired by geometry and mathematics and vice-versa. Common design motifs between wrought iron gates and Poncelet loci have stimulated us to look at the latter both from a geometric and an aesthetic perspective. As seen in \cref{fig:more-gates}, wrought iron {\em kraftwerk} is on a league of its own. An interesting question is if Poncelet loci could ever be used as a basis for new metalwork and/or architectural design. 

\subsubsection*{Acknowledgements}
\noindent I am grateful to my sister Regina Reznik for her artistic design work, and Iverton Darlan for being a co-developer of the locus app.

\appendix

\bibliographystyle{splncs04}
\bibliography{999_refs,999_refs_rgk}

\begin{thebibliography}{10}
\providecommand{\url}[1]{\texttt{#1}}
\providecommand{\urlprefix}{URL }
\providecommand{\doi}[1]{https://doi.org/#1}

\bibitem{dreher2015-art}
History of computer art. IASLonline NetArt (2015),
  \url{http://iasl.uni-muenchen.de/links/GCA_Indexe.html}

\bibitem{akopyan2020-invariants}
Akopyan, A., Schwartz, R., Tabachnikov, S.: Billiards in ellipses revisited.
  Eur. J. Math.  (9 2020), \url{doi:10.1007/s40879-020-00426-9}

\bibitem{bialy2020-invariants}
Bialy, M., Tabachnikov, S.: {Dan Reznik's} identities and more. Eur. J. Math.
  (9 2020), \url{doi:10.1007/s40879-020-00428-7}

\bibitem{bos1987}
Bos, H.J.M., Kers, C., Oort, F., Raven, D.W.: {P}oncelet's closure theorem.
  Exposition. Math.  \textbf{5}(4),  289--364 (1987)

\bibitem{centina15}
del Centina, A.: {P}oncelet's porism: a long story of renewed discoveries i.
  Arch. Hist. Exact Sci.  \textbf{70}(2),  1--122 (2016).
  \doi{10.1007/s00407-015-0163-y}, \url{doi.org/10.1007/s00407-015-0163-y}

\bibitem{caliz2020-area-product}
Chavez-Caliz, A.: More about areas and centers of {P}oncelet polygons. Arnold
  Math J.  (8 2020), \url{doi:10.1007/s40598-020-00154-8}

\bibitem{darlan2021-app}
Darlan, I., Reznik, D.: An app for visual exploration, discovery, and sharing
  of {P}oncelet 3-periodic phenomena. arXiv:2106.04521 (6 2021),
  \url{dan-reznik.github.io/ellipse-mounted-loci-p5js}

\bibitem{dragovic11}
Dragovi\'{c}, V., Radnovi\'{c}, M.: Poncelet Porisms and Beyond: Integrable
  Billiards, Hyperelliptic Jacobians and Pencils of Quadrics. Frontiers in
  Mathematics, Springer, Basel (2011)

\bibitem{corentin2021-circum}
Fierobe, C.: On the circumcenters of triangular orbits in elliptic billiard.
  Journal of Dynamical and Control Systems  \textbf{27}(4),  693--705 (2021).
  \doi{10.1007/s10883-021-09537-2}

\bibitem{garcia2019-incenter}
Garcia, R.: Elliptic billiards and ellipses associated to the 3-periodic
  orbits. American Mathematical Monthly  \textbf{126}(06),  491--504 (2019)

\bibitem{garcia2020-similarity-I}
Garcia, R., Reznik, D.: Related by similarity {I}: Poristic triangles and
  3-periodics in the elliptic billiard (06 2021)

\bibitem{garcia2020-family-ties}
Garcia, R., Reznik, D.: Family ties: Relating {P}oncelet 3-periodics by their
  properties. J. Croatian Soc. for Geom. Gr. {(KoG)}  (2022), to appear

\bibitem{garcia2020-new-properties}
Garcia, R., Reznik, D., Koiller, J.: New properties of triangular orbits in
  elliptic billiards. Amer. Math. Monthly  \textbf{128}(10),  898--910 (2021).
  \doi{10.1080/00029890.2021.1982360}

\bibitem{stachel2019-conics}
Glaeser, G., Stachel, H., Odehnal, B.: The Universe of Conics: From the ancient
  Greeks to 21st century developments. Springer (2016)

\bibitem{gutruf2010-vermeer}
Gutruf, G., Stachel, H.: The hidden geometry in {V}ermeer’s ``the art of
  painting''. J. for Geom. and Graphics  \textbf{14}(2),  187--202 (2010)

\bibitem{helman2021-power-loci}
Helman, M., Laurain, D., Garcia, R., Reznik, D.: Invariant center power and
  loci of {P}oncelet triangles. J. Dyn. \& Contr. Sys.  (10 2021),
  \url{doi:10.1007/s10883-021-09580-z}, to appear

\bibitem{helman2021-theory}
Helman, M., Laurain, D., Reznik, D., Garcia, R.: Poncelet triangles: a theory
  for locus ellipticity. Beitr. Algebra Geom.  (2022), to appear

\bibitem{karger98-classical}
Karger, A.: Classical geometry and computers. J. for Geom. and Graphics
  \textbf{2}(1),  7--15 (1998)

\bibitem{etc}
Kimberling, C.: Encyclopedia of triangle centers (2019), \url{bit.ly/3mOOver}

\bibitem{lordick2021-vermeer}
Lordick, D.: Parametric reconstruction of the space in {V}ermeer's painting
  ``{Girl Reading a Letter at an Open Window}''. J. for Geom. and Graphics
  \textbf{16}(1),  69--79 (2021)

\bibitem{odehnal2011-poristic}
Odehnal, B.: Poristic loci of triangle centers. J. Geom. Graph.
  \textbf{15}(1),  45--67 (2011)

\bibitem{pamfilos2004}
Pamfilos, P.: On some actions of {$D_3$} on a triangle. Forum Geometricorum
  \textbf{4},  157--176 (2004),
  \url{forumgeom.fau.edu/FG2004volume4/FG200420.pdf}

\bibitem{popper2003}
Popper, F.: Kinetic art. Grove art online (2003).
  \doi{10.1093/gao/9781884446054.article.T046632}, \url{https://bit.ly/3Aa2YJ0}

\bibitem{reznik2021-artful-long}
Reznik, D.: Artful loci of poncelet triangles (long version). {Google Slides}
  (2021), \url{bit.ly/3dSr4wP}

\bibitem{reznik2021-artful-short}
Reznik, D.: Artful loci of poncelet triangles (short version). {Google Slides}
  (2021), \url{bit.ly/34sO8Px}

\bibitem{reznik2021-artful-youtube}
Reznik, D.: Artful loci of poncelet triangles (video). {YouTube} (2021),
  \url{youtu.be/l-O5UT8tpuw}

\bibitem{reznik2021-locus-app-tutorial}
Reznik, D.: Locus app tutorial: a video walkthrough. {YouTube} (2021),
  \url{bit.ly/3iCKHxn}

\bibitem{reznik2020-similarityII}
Reznik, D., Garcia, R.: Related by similarity {II}: {P}oncelet 3-periodics in
  the homothetic pair and the {B}rocard porism. Intl. J. of Geom.
  \textbf{10}(4),  18--31 (03 2021)

\bibitem{reznik2020-ballet}
Reznik, D., Garcia, R., Koiller, J.: The ballet of triangle centers on the
  elliptic billiard. Journal for Geometry and Graphics  \textbf{24}(1),
  079--101 (2020)

\bibitem{reznik2020-intelligencer}
Reznik, D., Garcia, R., Koiller, J.: Can the elliptic billiard still surprise
  us? Math Intelligencer  \textbf{42},  6--17 (2020).
  \doi{10.1007/s00283-019-09951-2}

\bibitem{regina2021-private}
Reznik, R.: Private communication (2021)

\bibitem{olga14}
Romaskevich, O.: On the incenters of triangular orbits on elliptic billiards.
  Enseign. Math.  \textbf{60},  247--255 (2014). \doi{10.4171/LEM/60-3/4-2}

\bibitem{sharp2015-artzt}
Sharp, J.: Artzt parabolas of a triangle. The Mathematical Gazette
  \textbf{99}(546),  444--463 (2015)

\bibitem{skutin2013}
Skutin, A.: On rotation of a isogonal point. J. of Classical Geom.  \textbf{2},
   66--67 (2013)

\bibitem{stachel2021-billiards}
Stachel, H.: The geometry of billiards in ellipses and their {P}oncelet grids
  (2021)

\bibitem{mw}
Weisstein, E.: Mathworld. MathWorld--A Wolfram Web Resource  (2019),
  \url{mathworld.wolfram.com}

\bibitem{zaslavsky2001-poncelet}
Zaslavsky, A., Chelnokov, G.: The {P}oncelet theorem in euclidean and algebraic
  geometry. Mathematical Education  \textbf{4}(19),  49--64 (2001), in Russian

\bibitem{zaslavsky2003-trajectories}
Zaslavsky, A., Kosov, D., Muzafarov, M.: Trajectories of remarkable points of
  the {P}oncelet triangle. Kvant  \textbf{2},  22--25 (2003)

\end{thebibliography}

\end{document}